\documentclass[12pt]{amsart}
\usepackage{amssymb, amscd}

\begin{document}

\title[\null]
{Hodge metrics and the curvature of \\ higher direct images}

\author[\null]{Christophe Mourougane and Shigeharu TAKAYAMA} 

\begin{abstract} 
Using the harmonic theory developed by Takegoshi for representation of 
relative cohomology and the framework of computation of curvature of 
direct image bundles by Berndtsson, 
we prove that the higher direct images by a smooth morphism of 
the relative canonical bundle twisted by a semi-positive vector bundle 
are locally free and semi-positively curved, when endowed 
with a suitable Hodge type metric. 

\vskip5pt

\noindent
R\'esum\'e.
Nous utilisons la th\'eorie de repr\'esentation par formes harmoniques
des classes de cohomologie relative d\'evelopp\'ee par Takegoshi
et la structure des calculs de courbure de fibr\'es images directes
d\'evelopp\'ee par Berndtsson, pour \'etudier les images directes 
sup\'erieures
par un morphisme lisse du fibr\'e canonique relatif
tensoris\'e par un fibr\'e vectoriel holomorphe hermitien semi-positif.
Nous montrons qu'elles sont localement libres
et que, munies de m\'etriques convenables de type Hodge,
elles sont \`a courbure semi-positive.
\end{abstract}

\maketitle


\baselineskip=18pt


\theoremstyle{plain}
  \newtheorem{thm}{Theorem}[section]
  \newtheorem{main}[thm]{Main Theorem}
  \newtheorem{defthm}[thm]{Definition-Theorem}
  \newtheorem{prop}[thm]{Proposition}
  \newtheorem{lem}[thm]{Lemma}
  \newtheorem{cor}[thm]{Corollary}
  \newtheorem{conj}[thm]{Conjecture}
  \newtheorem{sublem}[thm]{Sublemma}
  \newtheorem{mainlem}[thm]{Main Lemma}
\theoremstyle{definition}
  \newtheorem{dfn}[thm]{Definition}
  \newtheorem{exmp}[thm]{Example}
  \newtheorem{co-exmp}[thm]{Counter-Example}
  \newtheorem{prob}[thm]{Problem}
  \newtheorem{notation}[thm]{Notation}
  \newtheorem{quest}[thm]{Question}
\theoremstyle{remark}
  \newtheorem{rem}[thm]{Remark}
  \newtheorem{com}[thm]{Comment}
\renewcommand{\theequation}{\thesection.\arabic{equation}}
\setcounter{equation}{0}

\newcommand{\BC}{{\mathbb{C}}}
\newcommand{\BN}{{\mathbb{N}}}
\newcommand{\BP}{{\mathbb{P}}}
\newcommand{\BQ}{{\mathbb{Q}}}
\newcommand{\BR}{{\mathbb{R}}}

\newcommand{\CalC}{{\mathcal{C}}}
\newcommand{\CalD}{{\mathcal{D}}}
\newcommand{\CE}{{\mathcal{E}}}
\newcommand{\CF}{{\mathcal{F}}}
\newcommand{\CH}{{\mathcal{H}}}
\newcommand{\CI}{{\mathcal{I}}}
\newcommand{\CJ}{{\mathcal{J}}}
\newcommand{\CK}{{\mathcal{K}}}
\newcommand{\CL}{{\mathcal{L}}}
\newcommand{\CM}{{\mathcal{M}}}
\newcommand{\CO}{{\mathcal{O}}}
\newcommand{\CS}{{\mathcal{S}}}
\newcommand{\CU}{{\mathcal{U}}}

\newcommand\ga{\alpha}
\newcommand\gb{\beta}
\def\th{\theta}
\newcommand\vth{\vartheta}
\newcommand\Th{\Theta}
\newcommand\ep{\varepsilon}
\newcommand\Ga{\Gamma}
\newcommand\w{\omega}
\newcommand\Om{\Omega}
\newcommand\Dl{\Delta}
\newcommand\del{\delta}
\newcommand\sg{\sigma}
\newcommand\vph{\varphi}
\newcommand\lam{\lambda}
\newcommand\Lam{\Lambda}

\newcommand\lra{\longrightarrow}
\newcommand\ra{\rightarrow}
\newcommand\ot{\otimes}
\newcommand\wed{\wedge}
\newcommand\ol{\overline}
\newcommand\isom{\, \wtil{\lra}\, }

\newcommand\ai{\sqrt{-1}}
\newcommand\rd{{\partial}}
\newcommand\rdb{{\overline{\partial}}}
\newcommand\levi{\ai\rd\rdb}

\newcommand\wtil{\widetilde}
\newcommand\what{\widehat}
\newcommand\hsg{\widehat \sg}

\newcommand\End{\mbox{{\rm End}}\, }
\newcommand\im{\mbox{{\rm im}}\, }
\newcommand\Ima{\mbox{{\rm Im}}\, }
\newcommand\Ker{\mbox{{\rm Ker}}\, }

\newcommand\Rq{R^qf_*\Om_{X/Y}^n(E)}
\newcommand\RqX{R^qf_*\Om_X^{n+m}(E)}


\section{Introduction}

This is a continuation of our works \cite{M}\,\cite{MT} on the metric
positivity of direct image sheaves of adjoint bundles.
The goal of this paper is to prove the following

\begin{thm} \label{MT}
Let $f : X \lra Y$ be a holomorphic map of complex manifolds,
which is smooth, proper, K\"ahler, surjective, and with connected fibers.
Let $(E, h)$ be a holomorphic vector bundle on $X$ with
a Hermitian metric $h$ of semi-positive curvature in the sense of Nakano.
Then for any $q \ge 0$, the direct image sheaf $R^qf_*\Om_{X/Y}^n(E)$ 
is locally free and Nakano semi-positive,
where $n$ is the dimension of fibers.
\end{thm}

A real $(1,1)$-form $\w$ on $X$ is said to be {\it a relative K\"ahler form}
for $f : X \lra Y$, if for every point $y \in Y$, 
there exists a local coordinate $(W; (t_1, \ldots, t_m))$ around $y$ 
such that $\w + cf^*(\ai\sum_j dt_j \wed d\ol{t_j})$ is a K\"ahler form
on $f^{-1}(W)$ for a constant $c$.
A morphism $f$ is said to be {\it K\"ahler}, if there exists 
a relative K\"ahler form $\w_f$ for $f$ (see \cite[6.1]{Tk}). 

In case when $E$ is a semi-positive line bundle and $q=0$, 
Theorem \ref{MT} is a theorem of Berndtsson \cite[1.2]{B}.
In our previous paper \cite{MT}, we obtained independently from \cite{B},
a weaker semi-positivity:\ the Griffiths semi-positivity
of $f_*\Om_{X/Y}^n(E)$ for a semi-ample line bundle $E$.
Right after two papers \cite{B}\,\cite{MT}, especially \cite{B} have appeared,
the analogous statement for higher direct images has been
considered as a next problem among others.
Theorem \ref{MT} solves this problem for Nakano semi-positive 
vector bundles $E$.

For a Nakano semi-positive vector bundle $E$,
the local freeness $R^qf_*\Om_{X/Y}^n(E)$ is a consequence
of Takegoshi's injectivity theorem \cite{Tk}.
Here is one point where we use the smoothness of $f$.
We can only expect the torsion freeness in general,
by Koll\'ar \cite{Ko86} (\cite{Tk} in analytic setting).
Another theorem in \cite{Tk} shows that 
$R^qf_*\Om_{X/Y}^n(E)$ can be embedded into $f_*\Om_{X/Y}^{n-q}(E)$ 
at least locally on $Y$, and that $R^qf_*\Om_{X/Y}^n(E) = 0$ for $q > n$.
The sheaf $f_*\Om_{X/Y}^{n-q}(E)$ has a natural Hermitian metric
induced from a relative K\"ahler form $\w_f$ and $h$
(at least on which it is locally free) the so-called Hodge metric.
For local sections $\sg, \tau \in H^0(W, f_*\Om_{X/Y}^{n-q}(E))$,
the inner product at $y \in W \subset Y$ is given by
$\int_{X_y} (\sg|_{X_y}) \wed \ol{*}_{h_y} (\tau|_{X_y})$,
where $X_y = f^{-1}(y)$ is the fiber, and $\ol{*}_{h_y}$
is the ``star''-operator with respect to 
$\w_y = \w_f|_{X_y}$ and $h_y = h|_{X_y}$.
By pulling back this Hodge metric 
via the injection, say
$S_f^q : R^qf_*\Om_{X/Y}^n(E) \lra f_*\Om_{X/Y}^{n-q}(E)$,
we have a Hermitian metric on $R^qf_*\Om_{X/Y}^n(E)$ in the usual sense.
Our original contribution is the right definition of this
``Hodge metric'' on $R^qf_*\Om_{X/Y}^n(E)$, and 
our main theorem is that its curvature is Nakano semi-positive.

The space $H^q(X_y, \Om_{X_y}^n(E_y))$ has another natural inner product 
with respect to $\w_y$ and $h_y$.
For cohomology classes $u_y, v_y \in H^q(X_y, \Om_{X_y}^n(E_y))$, 
it is given by 
$\int_{X_y} u_y' \wed \ol *_{h_y} v_y'$,
where $u_y'$ and $v_y'$ are the harmonic representatives
of $u_y$ and $v_y$ respectively.
These fiberwise inner products also define a Hermitian metric on 
$R^qf_*\Om_{X/Y}^n(E)$.
We first tried to compute its curvature, but we did not succeed it.

We follow \cite{B}, not \cite{MT}, for the method of computation
of the curvature.
Since one can directly see the original method in \cite{B},
let us explain how different from \cite{B}, i.e.,
the differences in cases $q = 0$ and $q > 0$.
In case $q = 0$, the map 
$S_f^0 : f_*\Om_{X/Y}^n(E) \lra f_*\Om_{X/Y}^{n}(E)$
is an isomorphism, in fact the multiplication by a constant.
Moreover $f_*\Om_{X/Y}^{n}(E)$ is locally free thanks to 
Ohsawa-Takegoshi type $L^2$-extension theorem \cite{OT}\,\cite{O}\,\cite{Ma}, 
and the $(1,0)$-derivative of $\sg \in H^0(Y, f_*\Om_{X/Y}^{n}(E))
= H^0(X, \Om_{X/Y}^{n}(E))$ vanishes on each fiber $X_y$
by simply a bidegree reason.
However in case $q > 0$, we have no local freeness of $f_*\Om_{X/Y}^{n-q}(E)$,
nor the vanishing of the $(1,0)$-derivative of 
$\sg \in H^0(Y, f_*\Om_{X/Y}^{n-q}(E)) = H^0(X, \Om_{X/Y}^{n-q}(E))$
on $X_y$.
To overcome these difficulties, we need to restrict ourselves to consider
the image of $S_f^q : R^qf_*\Om_{X/Y}^n(E) \lra f_*\Om_{X/Y}^{n-q}(E)$.
Then we have a local freeness as we mentioned above, and the vanishing of 
the $(1,0)$-derivative thanks to an estimate by Takegoshi \cite{Tk}.
This is a key fact, and which is like the $d$-closedness of
holomorphic $p$-forms on a compact K\"ahler manifold.
After getting those key observations:\ the local freeness,
the right Hodge metric to be considered, and the closedness of
holomorphic sections, the computation of the curvature 
is a straight forward generalization of \cite{B}.

There are many positivity results of direct image sheaves of 
relative canonical bundles and of adjoint bundles,
which are mostly about the positivity in algebraic geometry.
We will recall only a few here.
The origin is due to Griffiths in his theory on the variation
of Hodge structures \cite{Gr}.
Griffiths' work has been generalized by Fujita \cite{Ft}, Kawamata \cite{Ka}, 
Viehweg \cite{Vi1}, Koll\'ar \cite{Ko86},
and so on, in more algebro-geometric setting.
There are also positivity results on higher direct images by
Moriwaki \cite{Mw}, Fujino \cite{Fn} and so on.
We refer to \cite{Mori}\,\cite{EV}\,\cite{Vi2}\,\cite{Ko2} for further remarks
on related results.
On the analytic side, it can be understood as the plurisubharmonic variation
of related functions to the Robin constant \cite{Y}\,\cite{LY}, 
or of the Bergman kernels \cite{MY}.
There is also a series of works by Yamaguchi.
As he mentioned in \cite{B}, his method is inspired by those works.
There are also related recent works by Berndtsson-P\u{a}un \cite{BP}
and Tsuji \cite{Ts}.

{\it Acknowledgments}.
We would like to thank Professor Berndtsson for his correspondences 
in many occasions, for answering questions, and showing
us a revised version of his paper \cite{B}.
A part of this work was done during the second named author's stay in Rennes.
He would like to thank the mathematical department of Rennes
for a support to stay there.


\section{Preliminaries}

\subsection{Hermitian vector bundles} \label{hvb}

Let $X$ be a complex manifold of dimension $n$ with a Hermitian
metric $\w$, and let $E$ be a holomorphic vector bundle of rank $r$
on $X$ with a Hermitian metric $h$.
Let $(E^*, h^*)$ be the dual vector bundle.
Let $A^{p,q}(X, E)$ be the space of $E$-valued smooth $(p,q)$-forms,
and $A^{p,q}_0(X, E)$ be the space of $E$-valued smooth $(p,q)$-forms
with compact support.    
Let $* : A^{p,q}(X, E) \lra A^{n-q,n-p}(X, E)$ be the Hodge star-operator
with respect to $\w$.
For any $u \in A^{p,q}(X, E)$ and $v \in A^{s,t}(X, E)$,
we define $u \wed h \ol v \in A^{p+s,q+t}(X, \BC)$ as follows.
We take a local trivialization of $E$ on an open subset $U \subset X$,
and we regard $u = {}^t(u_1, \ldots, u_r)$ as a row vector
with $(p,q)$-forms $u_j$ on $U$, and similarly for $v = {}^t(v_1, \ldots, v_r)$.
The Hermitian metric $h$ is then a matrix valued function 
$h = (h_{jk})$ on $U$.
We define $u \wed h \ol v$ locally on $U$ by 
$$
	u \wed h\ol v=  \sum_{j,k} u_j \wed h_{jk} \ol{v_k} 
	\in A^{p+s,q+t}(U, \BC).
$$
We should write ${}^tu \wed h\ol v$, but if there is no risk of
confusions, we will write in this way.
In this manner, we can define anti-linear isomorphisms
$\sharp_h : A^{p,q}(X, E) \lra A^{q,p}(X, E^*)$ by
$\sharp_h u = h\ol u$, and 
$\ol*_h = \sharp_h \circ * : A^{p,q}(X, E) \lra A^{n-p,n-q}(X, E^*)$ by
$\ol*_h u = h\ol{*u}$.
The inner product on $A^{p,q}_0(X, E)$ is defined by
$(u, v)_h = \int_X u \wed \ol*_h v$.
Denote by $D_h = \rd_h + \rdb$ the metric connection,
and by $\Th_h = D_h^2$ the curvature of $(E, h)$.
The Hermitian vector bundle $(E, h)$ is said to be 
{\it Nakano semi-positive} (resp.\ {\it Nakano positive}),
if the $\End (E)$-valued real $(1,1)$-from $\ai\Th_h$
is positive semi-definite (resp.\ positive definite) quadratic form 
on each fiber of the vector bundle $T_X \ot E$.

We define $\vth_h : A^{p,q}(X, E) \lra A^{p,q-1}(X, E)$ by 
$\vth_h = - * \rd_h * = - \ol *_{h^*} \rdb \ol *_h$, which is
the formal adjoint operator of $\rdb : A^{p,q}(X, E) \lra A^{p,q+1}(X, E)$
with respect to the inner product $(\ ,\ )_h$.
We also define $\ol\vth : A^{p,q}(X, E) \lra A^{p-1,q}(X, E)$ by 
$\ol\vth = -* \rdb *$, which is the formal adjoint operator of 
$\rd_h : A^{p,q}(X, E) \lra A^{p+1,q}(X, E)$
with respect to the inner product $(\ ,\ )_h$.
We denote by $e(\th)$ the left exterior product acting on
$A^{p,q}(X, E)$ by a form $\th \in A^{s,t}(X, \BC)$.
Then the adjoint operator $e(\th)^*$ with respect to 
the inner product $(\ ,\ )_h$ is defined by
$e(\th)^* = (-1)^{(p+q)(s+t+1)}* e(\ol{\th})*$.
For instance we set $\Lam_\w = e(\w)^*$.
We recall the following very useful relation
(\cite[1.2.31]{Huy} \cite[6.29]{Vo}):\ 

\begin{lem} \label{primitive}
For a primitive element $u \in A^{p,q}(X, E)$, 
i.e., $p+q=k \leq n$ and $\Lam_\w u = 0$, 
the Hodge $*$-operator reads
$$
	* (\w^j \wed u) 
	= \ai^{p-q}(-1)^{\frac{k(k+1)}{2}} \frac{j!}{(n-k-j)!} 
		\w^{n-k-j} \wed u
$$
for every $0 \le j \le n-k$.
\end{lem}

As immediate consequences, we have

\begin{cor} \label{relation}
Denote by $c_{n-q} = \ai^{(n-q)^2} = \ai^{n-q} (-1)^{(n-q)(n-q-1)/2}$.

(1) Let $a, b \in A^{n-q, 0}(X, E)$.
Then $*a = (\ol{c_{n-q}}/q!) \w^q \wed a$,
$\ol{*}_h a = (c_{n-q}/q!) \w^q \wed h \ol a$,
$a = (c_{n-q}/q!) *(\w^q \wed a)$,
and 
$a \wed \ol{*}_h b = (c_{n-q}/ q!) \w^q \wed a \wed h \ol{b}$.

(2) Let $u \in A^{n,q}(X, E)$.
Then $u = (c_{n-q}/q!) \w^q \wed * u$.
\end{cor}

\begin{notation}
We use the following conventions often.
Denote by $c_d = \ai^{d^2} = \ai^d (-1)^{d(d-1)/2}$ 
for any non-negative integer $d$.
Let $t = (t_1, \ldots, t_m)$ be the coordinates of $\BC^m$.

(1)
$dt = dt_1 \wed \ldots \wed dt_m$, 
$d\ol{t} = \ol{dt}$, and
$dV_t := c_m dt \wed d\ol{t} 
= \bigwedge_{j=1}^m \ai dt_j \wed d\ol{t_j} > 0$.

(2)
Let $\what{dt_j}$ be a smooth $(m-1, 0)$-form 
without $dt_j$ such that $dt_j \wed \what{dt_j} = dt$, and 
$\what{d\ol{t_j}} = \ol{\what{dt_j}}$.

(3)
Let $\what{dt_j \wed d\ol{t_k}}$ be a smooth $(m-1, m-1)$-form 
without $dt_j$ and $d\ol{t_k}$ such that
$\ai dt_j \wed d\ol{t_k} \wed \what{dt_j \wed d\ol{t_k}} = c_m dt \wed d\ol t$.
\end{notation}

\subsection{Set up} \label{setup}

In the rest of this paper, we will use the following set up.

Let $X$ and $Y$ be complex manifolds of $\dim X = n+m$
and $\dim Y = m$.
Let $f : X \lra Y$ be a holomorphic map, which is 
smooth, proper, K\"ahler, surjective, and with connected fibers.
Let $(E, h)$ be a holomorphic vector bundle on $X$ of rank $r$, 
with a Hermitian metric $h$ whose curvature $\Th_h$ is semi-positive
in the sense of Nakano.

(I) a general setting:\ $f : (X, \w_f) \lra Y$.
We take a relative K\"ahler form $\w_f$ for $f$,
and let $\kappa_f = \{\w_f\}$ be the de Rham cohomology class.
On each fiber $X_y$, we have a K\"ahler form $\w_y = \w_f|_{X_y}$,
and a Nakano semi-positive vector bundle $(E_y, h_y) = (E, h)|_{X_y}$.


(II) a localized setting of (I):\ $f : (X, \w) \lra Y \subset \BC^m$.
We further assume that the base $Y$ is a unit ball in $\BC^m$ 
with coordinates $t = (t_1, \ldots, t_m)$ and 
with admissible charts over $Y$ (see below).
We take a global K\"ahler form 
$\w = \w_f + cf^*(\ai \sum dt_j \wed d\ol{t_j})$ on $X$ for large $c > 0$,
without changing the class $\kappa_f$ 
nor the fiberwise K\"ahler forms $\w_y$.

Since $f : X \lra Y$ is smooth, for every point $y \in Y$, 
we can take a local coordinate
$(W; t = (t_1, \ldots, t_m))$ centered at $y$, so that
$(W; t)$ is a unit ball in $\BC^m$, and a system of local coordinates
$\CU = \{(U_\ga; z_\ga, t); \ \ga = 1,2,3, \ldots \}$ of $f^{-1}(W)$ 
which is locally finite, and every $U_\ga$ is biholomorphic to 
a product $D_\ga \times W$ for a domain $D_\ga$ in $\BC^n$;
$x \mapsto (z_\ga^1(x), \ldots, z_\ga^n(x), t_1, \ldots, t_m)$, 
namely the projection from $D_\ga \times W$ to $W$ is compatible 
with the map $f|_{U_\ga}$.
We can write
$z_\ga^j = f_{\ga\gb}^j(z_\gb^1, \ldots, z_\gb^n, t_1, \ldots, t_m)$
for $1 \le j \le n$ on $U_\ga \cap U_\gb$.
All $f_{\ga\gb}^j(z_\gb, t)$ are holomorphic in $z_\gb$ and $t$.
We call it {\it admissible charts} $\CU$ over $W$ 
(cf.\ \cite[\S 2.3]{Kd}).

Since our assertions are basically local on $Y$,
we will mostly use the set up (II).
The set up (I) will be used in subsections 3.1, 3.3, 4.3 and \ref{defhod}.


\section{Generalities of relative differential forms}

Let $f : (X, \w_f) \lra Y$ and $(E, h)$ be as in \S \ref{setup}.I.
We recall the complex analytic properties of
the relative cotangent bundle $\Om_{X/Y} = \Om_X/f^*\Om_Y$ 
and the bundle of relative holomorphic $p$-forms 
$\Om_{X/Y}^p = \bigwedge^p \Om_{X/Y}$.
We will not distinguish a vector bundle and 
the corresponding locally free sheaf.
For a subset $S \subset Y$, we denote by $X_S = f^{-1}(S)$ 
and $E_S = E|_{X_S}$.

\subsection{Definition of relative differential forms}

Let $U \subset X$ be an open subset. 

(1)
For a form $u \in A^{p,q}(U,E)$, we have the restriction 
$u|_{X_y \cap U} \in A^{p,q}(X_y \cap U, E)$ on each 
fiber over $y \in Y$,
which is the pull-back as a form via the inclusion $X_y \lra X$.
Two forms $u, v \in A^{p,q}(U,E)$ are said to be {\it $f$-equivalent} 
``$u \sim v$'', if
$u|_{X_y \cap U} = v|_{X_y \cap U}$ for any $t \in Y$.
We denote the set of equivalence classes by
$$
	A^{p,q}(U/Y, E) = A^{p,q}(U,E)/\sim.
$$
The set $A^{p,q}(U/Y, E)$ will be called the space of relative 
differential forms on $U$.
We denote by $[u] \in A^{p,q}(U/Y, E)$ the equivalence class of
$u \in A^{p,q}(U, E)$.

(2)
A form $u \in A^{p,0}(U, E)$ is said to be
{\it holomorphic on each fiber}, if 
the restriction $u|_{X_y}$ is holomorphic, i.e., 
$u|_{X_y} \in H^0(X_y, \Om_{X_y}^p(E_y))$ for every $y \in Y$. 
A form $u \in A^{p,0}(U, E)$ is said to be {\it relatively holomorphic}, 
if for any local chart $(W; t=(t_1, \ldots, t_m))$ of $Y$,
the form $u \wed f^*dt$ is holomorphic on $X_W \cap U$, i.e.,
$u \wed f^*dt \in H^0(X_W \cap U, \Om_X^{p+m}(E))$.

%

(3)
For a function $\ga \in A^0(Y, \BC)$ and $[u] \in A^{p,0}(U/Y, E)$, 
we can define $\ga[u] := [(f^*\ga)u] \in A^{p,0}(U/Y, E)$.
For each open subset $W \subset Y$, we set
\begin{equation*} 
\begin{aligned}
A^0(W, f_*\Om_{X/Y}^p(E)) & :=
\{ [u] \in A^{p,0}(X_W/W, E);\ u \text{ is holomorphic on each fiber} \}, \\
H^0(W, f_*\Om_{X/Y}^p(E)) & := 
\{ [u] \in A^{p,0}(X_W/W, E);\ u \text{ is relatively holomorphic} \}.
\end{aligned}
\end{equation*} 
We can see that 
$A^0(W, f_*\Om_{X/Y}^p(E))$ becomes an $A^0(W, \BC)$-module, and
$H^0(W, f_*\Om_{X/Y}^p(E))$ becomes an $H^0(W, \CO_W)$-module.

(4)
It is some times convenient to use local coordinates
to look at those properties above.
Let $u \in A^{p,q}(X, E)$.
On an admissible chart $(U_\ga; z_\ga, t)$ as above, we can write 
$$
	u = \sum_{I \in I_p, J \in I_q} 
		u_{IJ\ga} dz_\ga^I \wed d\ol{z_\ga^J} + R,
$$
where $u_{IJ\ga} = u_{IJ\ga}(z_\ga, t) \in A^0(U_\ga, \BC^r)$,
and $R \in \sum_j A^{p-1,q}(U_\ga, \BC^r) \wed dt_j 
+ \sum_j A^{p,q-1}(U_\ga, \BC^r) \wed d\ol{t_j}$. 
Here we use a standard convention.
We set $I_p = \{\{i_1, i_2, \ldots, i_p\}; \ 
1 \le i_1 < i_2 < \ldots < i_p \le n \}$, and $I_0$ is empty. 
For $I = \{i_1, i_2, \ldots, i_p\} \in I_p$, we denote by
$dz_\ga^I = dz_\ga^{i_1} \wed \ldots \wed dz_\ga^{i_p}$.
Similar for $J \in I_q$ and $d\ol{z_\ga^J}$.
%
%
%
The restriction on a fiber is locally given by 
$$
	u|_{X_y} = \sum_{I \in I_p, J \in I_q} 
		  u_{IJ\ga}|_{X_y}  dz_\ga^I \wed d\ol{z_\ga^J}.
$$
In particular, for two forms $u, v \in A^{p,q}(X, E)$,
they are $f$-equivalent $u \sim v$ if and only if
$u_{IJ\ga} = v_{IJ\ga}$ for any $(I,J) \in I_p \times I_q$
on any admissible chart $(U_\ga; z_\ga,t)$.

(5)
Let $u \in A^{p,0}(X, E)$.
On an admissible chart $(U_\ga; z_\ga, t)$, we have 
$$
	u = \sum_{I \in I_p} u_{I\ga} dz_\ga^I + R
$$
as above.
Therefore
$u$ is holomorphic on each fiber (resp.\ relatively holomorphic)
if and only if every $u_{I\ga}$ is holomorphic in $z_\ga$ 
(resp.\ holomorphic in $z_\ga$ and $t$)
for any $I \in I_p$ on any admissible chart $(U_\ga; z_\ga,t)$.

\subsection{Holomorphic structure of $f_*\Om_{X/Y}^p(E)$}

We also give the holomorphic structure on $f_*\Om_{X/Y}^p(E)$
by an action of the $\rdb$-operator.
Let $(W; (t_1, \ldots, t_m)) \subset Y$ be a local chart,
over which we have admissible charts.
Let $[\sg] \in A^0(W, f_*\Om_{X/Y}^p(E))$, which can be seen as 
a differentiable family of holomorphic forms.
Since $(\rdb \sg)|_{X_y} = \rdb(\sg|_{X_y}) = 0$,
we can write as
$$
	\rdb\sg = \sum_j \eta^j \wed dt_j + \sum_j \nu^j \wed d\ol{t_j}
$$
with some $\eta^j \in A^{p-1,1}(X_W, E)$ and some $\nu^j \in A^{p,0}(X_W, E)$.
In particular $\rdb(\sg \wed dt) = \sum_j \nu^j \wed d\ol{t_j} \wed dt$.
On an admissible chart $(U; z,t) = (U_\ga; z_\ga,t)$, we write 
$\sg = \sum_{I \in I_{p}} \sg_I dz^I + R$ with 
$R \in \sum_j A^{p-1,0}(U_\ga, \BC^r) \wed dt_j$.
Then we have 
$\nu^j|_{X_y \cap U} 
= (-1)^p \sum_{I \in I_p} (\rd \sg_I/\rd\ol{t_j})|_{X_y \cap U} dz^I$
$\in H^0(X_y \cap U, \Om_{X_y}^{p}(E_y))$ for every $j$.
In particular, the class $[\nu^j] \in A^0(W, f_*\Om_{X/Y}^p(E))$ is well-defined for $[\sg]$.
For $[\sg] \in A^0(W, f_*\Om_{X/Y}^p(E))$, we define 
$$
	\rdb[\sg] = \sum_j [\nu^j]d\ol{t_j} \ \in A^{0,1}(W, f_*\Om_{X/Y}^p(E)). 
$$
Here $A^{0,1}(W, f_*\Om_{X/Y}^p(E)) 
= A^0(W, f_*\Om_{X/Y}^p(E)) \ot A^{0,1}(W, \BC)$ as $A^0(W, \BC)$-module,
but it has only a formal meaning.
Then, $[\sg] \in H^0(W, f_*\Om_{X/Y}^p(E))$ if and only if 
$\rdb[\sg] \equiv 0$.
In fact both of them are characterized by the holomorphicity of 
all $\sg_I$ in $z$ and $t$ locally.

\begin{lem} \label{rdb}
Let $(W; (t_1, \ldots, t_m)) \subset Y$ be a local chart as above 
in \S \ref{setup}.
Let $\sg \in A^{p,0}(X_W, E)$ such that $[\sg] \in H^0(W, f_*\Om_{X/Y}^p(E))$.
Then
(1)
$$
	\rdb\sg = \sum_j \eta^j \wed dt_j 
$$
with some $\eta^j \in A^{p-1,1}(X_W, E)$,

(2) these $\eta^j$ are not unique, but $[\eta^j] \in A^{p-1,1}(X_W/W, E)$
are well-defined for $\sg$, 

(3) all $\eta^j|_{X_y}$ are $\rdb$-closed on any $X_y$, and

(4) \cite[Lemma 4.3]{B}
all $\eta^j|_{X_y} \wed \w_y^{q+1}$ are $\rdb$-exact on any $X_y$.
\end{lem}

\begin{proof}
(1) is now clear. 
We show (2) and (3).
For each $j$, $\sg \wed \what{dt_j} \in A^{p+m-1,0}(X_W, E)$
is well-defined for $\sg$, and so is 
$\rdb(\sg \wed \what{dt_j}) = \eta^j \wed dt$. 
Hence $[\eta^j]$ are well-defined for $\sg$, by Remark \ref{eq} below.
Moreover $(\rdb \eta^j) \wed dt = \rdb\,\rdb(\sg \wed \what{dt_j}) = 0$.
Hence we obtain $\rdb(\eta^j|_{X_y}) = (\rdb \eta^j)|_{X_y} = 0$
by Remark \ref{eq} again.

(4)
We fix $j$.
By a bidegree reason, we can write as
$\sg \wed \what{dt_j} \wed \w^{q+1} = a^j \wed dt$
with some $a^j \in A^{n, q+1}(X, E)$.
We note that the class $[a^j] \in A^{n, q+1}(X/Y, E)$ is well-defined
by Remark \ref{eq}.
By taking $\rdb$, we have $\eta^j \wed dt \wed \w^{q+1} = (\rdb a^j) \wed dt$.
Then $[\eta^j \wed \w^{q+1}] = [\rdb a^j]$ in $A^{n, q+2}(X/Y, E)$
by Remark \ref{eq}, and hence 
$\eta^j|_{X_y} \wed \w_y^{q+1} = \rdb (a^j|_{X_y})$ on any $X_y$.
\end{proof}

\begin{rem} \label{eq}
For $u, v \in A^{p,q}(X_W, E)$, a relation
$u \wed dt = v \wed dt$ implies $[u] = [v]$ in $A^{p,q}(X_W/W, E)$,
and the converse holds true in case $q = 0$.
\end{rem}

\begin{rem}
Each cohomology class $\{\eta^j|_{X_y}\} \in H^{p-1,1}(X_y, E_y)$
is well-defined for $[\sg] \in H^0(W, f_*\Om_{X/Y}^p(E))$.
As it is well-known, $\{\eta^j|_{X_y}\}$ is obtained by
the cup-product (up to a sign) with the Kodaira-Spencer class
of $f : X \lra Y$ at $y \in Y$ for $t_j$-direction.
However we will not use these remarks.
\end{rem}

\subsection{Canonical pairing} \label{can}

There is a canonical pairing on each stalk $f_*\Om_{X/Y}^p(E)_y$
with respect to $\w_y$ and $h_y$, via the natural inclusion
$f_*\Om_{X/Y}^p(E)_y \subset H^0(X_y, \Om_{X_y}^{p}(E_y))$.
At each point $y \in Y$, we have the fiberwise inner product
$$
	g_y(\sg_y, \tau_y) := (\sg_y, \tau_y)_{h_y} 
	= \int_{X_y} (c_{p}/(n-p)!) \w_y^{n-p} \wed \sg_y \wed h\ol{\tau_y}
$$
for $\sg_y, \tau_y \in A^{p,0}(X_y, E_y)$. 
As germs $\sg_y, \tau_y \in f_*\Om_{X/Y}^p(E)_y$,
we will denote by $g_y(\sg_y, \tau_y)$.
On the other hand, as forms $\sg_y, \tau_y \in H^0(X_y, \Om_{X_y}^{p}(E_y))$,
we will denote by $(\sg_y, \tau_y)_{h_y}$.
These two are the same, but our standing points are different, i.e.,
at a point $y \in Y$, or on the fiber $X_y$.

For relative forms $[\sg], [\tau] \in A^{p,0}(X_W/W, E)$
(or $H^0(W, f_*\Om_{X/Y}^{p}(E))$) over an open subset $W \subset Y$, 
the above fiberwise inner product gives
$$
   g([\sg],[\tau]) := f_*((c_{p}/(n-p)!) \w_f^{n-p} \wed \sg \wed h\ol\tau),
$$
where the right hand side is a push-forward as a current.
For a test $(m,m)$-form $\vph$ on $W$
(i.e., a smooth form with compact support),
we have
$f_*((c_p/(n-p)!) \w_f^{n-p} \wed \sg \wed h\ol\tau)(\vph)
:= \int_X (c_p/(n-p)!) \w_f^{n-p} \wed \sg \wed h\ol\tau \wed f^*\vph$.
Hence the right hand side does not depend on representatives 
$\sg$ nor $\tau$ (see Remark \ref{eq}).
Since the map $f$ is smooth, 
$g([\sg],[\tau])$ is in fact a smooth function on $W$.

We remark that the definition of the pairing $g$ depends only on
the fiberwise K\"ahler forms $\{\w_y\}_{y \in Y}$.
For example, over a local chart $(W; t) \subset Y$, we can replace
$\w_f$ by $\w_f + c f^*(\ai\sum dt_j \wed d\ol{t_j})$ for any $c \in \BR$
in the definition of $g([\sg],[\tau])$.
The pairing $g$ defines a Hermitian metric on
every locally free subsheaf of $f_*\Om_{X/Y}^p(E)$ 
in the usual sense, which we call the Hodge metric.
As a matter of fact, $g$ itself is called the Hodge metric 
on $f_*\Om_{X/Y}^p(E)$ commonly,
although it may not be locally free.


\section{Harmonic theory for Nakano semi-positive vector bundles}

We collect some fundamental results of Takegoshi \cite{Tk},
and immediate consequences from them.

\subsection{Absolute setting} \label{abs}

Let $(X_0, \w_0)$ be an $n$-dimensional compact K\"ahler manifold
and let $(E_0, h_0)$ be a holomorphic Hermitian vector bundle on $X_0$. 
We have an inner product $(\ ,\ )_{h_0}$ and the associated 
norm $\|\ \|_{h_0}$ on each $A^{p,q}(X_0, E_0)$.
Let $\CH^{p,q}(X_0, E_0)$ be the space of harmonic $(p,q)$-forms.
%
From Bochner-Kodaira-Nakano formula, it then follows that 
if $(E_0, h_0)$ is furthermore 
Nakano semi-positive and Nakano positive at one point,
then $H^q(X_0, \Om_{X_0}^n(E_0))$ vanish for all $q>0$.
%

Enoki \cite{E} and Takegoshi \cite{Tk} (a special case of 
Theorem \ref{Tk4.3} below) 
show that if $(E_0, h_0)$ is Nakano semi-positive,
then the Hodge $*$-operator yields injective homomorphism
$*_0 : \CH^{n,q}(X_0, E_0) \lra H^0(X_0, \Om_{X_0}^{n-q}(E_0))$. 
Recalling that $(c_{n-q}/q!) \w_0^q \wed *_0 u = u$ for 
$u \in A^{n,q}(X_0,E_0)$, it then follows that the Lefschetz operator 
$L_0^q : H^0(X_0, \Om_{X_0}^{n-q}(E_0)) \lra 
H^q(X_0, \Om_{X_0}^n(E_0))$ is surjective.
Hence we have $H^0(X_0, \Om_{X_0}^{n-q}(E_0)) 
= \Ker L_0^q \oplus *_0 \CH^{n,q}(X_0,E_0)$.

\subsection{Localized relative setting}

Let $f : (X, \w) \lra Y \subset \BC^m$ and $(E, h)$ be as in \S \ref{setup}.II.
We take a $C^\infty$ plurisubharmonic exhaustion function
$\Phi = f^*\sum_{j=1}^m |t_j|^2$ on $X$.
We take any $0 \le q \le n$.
Following \cite[4.3]{Tk}, we set the following subspace of 
$E$-valued harmonic $(n+m,q)$-forms with respect to $\w$ and $h$:
$$
	\CH^{n+m,q}(X,E,\Phi) = \{ u \in A^{n+m,q}(X, E) ; \
	\rdb u = \vth_h u = 0 \text{ and } e(\rdb\Phi)^*u = 0 \text { on } X \}.
$$
By \cite[4.3.i]{Tk}, $u \in \CH^{n+m,q}(X, E, \Phi)$
if and only if $\ol\vth u = 0$, 
$(\ai e(\Th_h + \rd\rdb \Phi) \Lam_\w u) \wed h\ol{u} = 0$
and $e(\rdb\Phi)^*u = 0$ on $X$.
One can check easily that $(f^*\ga)u$ satisfies
those latter three conditions, if $\ga \in H^0(Y, \CO_Y)$ and 
if $u \in \CH^{n+m,q}(X, E, \Phi)$.

\begin{thm} \label{Tk4.3}
\cite[4.3]{Tk}.
(1) 
The space $\CH^{n+m,q}(X,E,\Phi)$ does not depend on
$C^\infty$ plurisubharmonic exhaustion functions $\Phi$,
and has a natural structure of $H^0(Y, \CO_Y)$-module.

(2)
For $u \in \CH^{n+m,q}(X,E,\Phi)$, one has
$\rdb * u = 0$ and $\rd_h * u = 0$. 
In particular, the Hodge $*$-operator yields an injective
homomorphism 
$* : \CH^{n+m,q}(X, E, \Phi) \lra H^0(X, \Om_X^{n+m-q}(E))$,
and $\CH^{n+m,q}(X, E, \Phi)$ becomes a torsion free 
$H^0(Y, \CO_Y)$-module.
\end{thm}

Let $\iota' : Z^{n+m,q}_\rdb(X, E) \lra H^{q}(X, \Om_X^{n+m}(E))$
be the quotient map which induces the Dolbeault's isomorphism. 

\begin{thm} \label{Tk5.2} \cite[5.2]{Tk}. 
(1)
The space $\CH^{n+m,q}(X, E, \Phi)$ represents $H^q(X, \Om_X^{n+m}(E))$ 
as a torsion free $H^0(Y, \CO_Y)$-module,
in particular the quotient map $\iota'$ induces an isomorphism
$\iota : \CH^{n+m,q}(X, E, \Phi) \lra H^{q}(X, \Om_X^{n+m}(E))$.

(2) 
The injective homomorphism 
$* : \CH^{n+m,q}(X, E, \Phi) \lra H^0(X, \Om_X^{n+m-q}(E))$
induces a splitting homomorphism (up to a constant)
$$
* \circ \iota^{-1} : 
	H^{q}(X, \Om_X^{n+m}(E)) \lra H^0(X, \Om_X^{n+m-q}(E))
$$
for the Lefschetz homomorphism
$$
L^q : H^0(X, \Om_X^{n+m-q}(E)) \lra H^q(X, \Om_X^{n+m}(E)).
$$
such that $(c_{n+m-q}/q!) L^q \circ * \circ \iota^{-1} = id$.

(3)
Let $u \in \CH^{n+m,q}(X, E, \Phi)$.
Then the form $*u \in H^0(X, \Om_X^{n+m-q}(E))$ 
is saturated in base variables, i.e., $*u =  \sg_u \wed dt$ 
for some  $[\sg_u] \in H^0(X, \Om_{X/Y}^{n-q}(E))$
(see the proof of \cite[5.2.ii]{Tk}).
In particular, $u = (c_{n+m-q}/q!) \w^q \wed \sg_u \wed dt$ 
and the map $u \mapsto [\sg_u]$ is well-defined.
Thus the Hodge $*$-operator induces 
an injective homomorphism
$$
	S^q : \CH^{n+m,q}(X, E, \Phi) 
		\lra H^0(X, \Om_{X/Y}^{n-q}(E)).
$$
\end{thm}

In Theorem \ref{Tk5.2}\,(3), we used our assumption that $f$ is smooth.

We take a trivialization $\CO_Y \isom \Om_Y^m$ given by $1 \mapsto dt$,
which induces isomorphisms of sheaves
$\Om_{X/Y}^{n} \cong \Om_{X/Y}^{n} \ot f^*\Om_Y^m
\cong \Om_X^{n+m}$ by $[u] \mapsto u \wed dt$, 
and hence of cohomology groups
$\ga^q : H^q(X, \Om_{X/Y}^{n}(E)) \isom H^q(X, \Om_X^{n+m}(E))$. 
We also have an injection $\Om_{X/Y}^{n-q} \lra \Om_X^{n+m-q}$ by
$[\sg] \mapsto \sg \wed dt$, and hence an injection
$\gb^0 : H^0(X, \Om_{X/Y}^{n-q}(E)) \lra H^0(X, \Om_X^{n+m-q}(E))$. 
Combining with Theorem \ref{Tk5.2}, we have 
\begin{equation*} 
\begin{aligned}
\iota^{-1} \circ \ga^q & : 
	H^q(X, \Om_{X/Y}^{n}(E)) \isom 
	H^q(X, \Om_X^{n+m}(E)) \isom \CH^{n+m,q}(X, E, \Phi), 
\\
* = \gb^0 \circ S^q & : 
	\CH^{n+m,q}(X, E, \Phi) \lra
	H^0(X, \Om_{X/Y}^{n-q}(E)) \lra H^0(X, \Om_X^{n+m-q}(E)), 
\\
(\ga^q)^{-1}\circ L^q & : 
	H^0(X, \Om_X^{n+m-q}(E)) \lra
	H^q(X, \Om_X^{n+m}(E)) \isom H^q(X, \Om_{X/Y}^{n}(E)).
\end{aligned}
\end{equation*} 
Then Theorem \ref{Tk5.2}\,(2) reads the following relative version:\

\begin{cor} \label{relTk5.2}
Let
\begin{equation*} 
\begin{aligned}
S_f^q = S^q \circ \iota^{-1} \circ \ga^q & : 
	H^q(X, \Om_{X/Y}^{n}(E)) \lra H^0(X, \Om_{X/Y}^{n-q}(E)),
\\
L_f^q = (\ga^q)^{-1} \circ L^q \circ \gb^0 & : 
	H^0(X, \Om_{X/Y}^{n-q}(E)) \lra H^q(X, \Om_{X/Y}^{n}(E)).
\end{aligned}
\end{equation*} 
Then $(c_{n+m-q}/q!) L_f^q \circ S_f^q = id$
on $H^q(X, \Om_{X/Y}^{n}(E))$.
\end{cor}

We can also see, thanks to \cite[5.2.iv]{Tk} (see also \cite[6.5.i]{Tk}) 
that those constructions can be localized on $Y$, 
and induce homomorphisms of direct image sheaves.

\begin{cor} \label{split}
There exist homomorphisms induced from the Hodge $*$-operator and 
the Lefschetz homomorphism:\
$$
S_f^q : R^qf_*\Om_{X/Y}^{n}(E) \lra f_*\Om_{X/Y}^{n-q}(E),
\ \ \ 
L_f^q : f_*\Om_{X/Y}^{n-q}(E) \lra R^qf_*\Om_{X/Y}^{n}(E)
$$
so that $(c_{n+m-q}/q!) L_f^q \circ S_f^q = id$ on $R^qf_*\Om_{X/Y}^{n}(E)$.
In particular
$$
	f_*\Om_{X/Y}^{n-q}(E) = F^{n-q} \oplus \CK^{n-q}, 
	\text{ with } F^{n-q} = \Ima S_f^q \text{ and } \CK^{n-q} = \Ker L_f^q.
$$
\end{cor}


We translate some results above into explicite forms.

\begin{lem} \label{rdh}
Let $\sg \in A^{n-q,0}(X, E)$ such that $[\sg] \in H^0(Y, F^{n-q})$.
Then (1)
$$
	\rd_h \sg = \sum_j \mu^j \wed dt_j
$$
for some $\mu^j \in A^{n-q,0}(X, E)$,

(2) these $\mu^j$ are not unique, but $[\mu^j] \in A^{n-q,0}(X/Y, E)$
are well-defined for $\sg$, and

(3) $\rd_{h_y}(\mu^j|_{X_y}) = 0$ on any $X_y$ and all $j$.
\end{lem}

\begin{proof}
There exists $u \in \CH^{n+m,q}(X,E,\Phi)$ such that
$*u = \sg \wed dt \in H^0(X, \Om_X^{n+m-q}(E))$.
By Takegoshi:\ Theorem \ref{Tk4.3}, we have $\rd_h * u = 0$.
Hence $(\rd_h \sg) \wed dt = \rd_h * u = 0$, and we have (1).
We can show (2) and (3) by the same method in Lemma \ref{rdb}.
\end{proof}

\begin{rem}
Unlike in the case $q=0$ that is treated by degree considerations
\cite[\S 4]{B}, we used the semi-positivity here.
In general, for $[\sg] \in H^0(Y, f_*\Om_{X/Y}^{n-q}(E))$ with $q > 0$,
we can not have 
$\rd_h \sg = \sum \mu^j \wed dt_j$ for some $\mu^j \in A^{n-q,0}(X, E)$.
This is in fact a key property, and it makes various computations possible.
We also note that $\eta^j$ in Lemma \ref{rdb} and $\mu^j$ are not
well-defined for a class $[\sg] \in H^0(Y, F^{n-q})$, and that means,
we have some freedom of choices in a class.
\end{rem}

Stalks or fibers at a point $y \in Y$ will be denoted by
$f_*\Om_{X/Y}^{n-q}(E)_y, F^{n-q}_y, \CK^{n-q}_y$ respectively.
Those stalks can be seen as subspaces of $H^0(X_y, \Om_{X_y}^{n-q}(E_y))$,
i.e., $F^{n-q}_y \oplus \CK^{n-q}_y = f_*\Om_{X/Y}^{n-q}(E)_y 
\subset H^0(X_y, \Om_{X_y}^{n-q}(E_y))$.


\begin{lem} \label{ortho}
Let $\sg_y \in F^{n-q}_y$ and $\tau_y \in \CK^{n-q}_y$,
and regard them as elements of $H^0(X_y, \Om_{X_y}^{n-q}(E_y))$.
Then,
(1) $\rd_{h_y} \sg_y = 0$ in $A^{n-q+1,0}(X_y, E_y)$,

(2) $\w_y^q \wed \tau_y \in A^{n,q}(X_y, E_y)$ is $\rdb$-exact, and

(3) $(\sg_y, \tau_y)_{h_y} 
= \int_{X_y} (c_{n-q}/q!) \w_y^q \wed \sg_y \wed h_y \ol{\tau_y} = 0$.
\end{lem}

\begin{proof}
We will argue at $y=0$.

(1)
Since $Y$ is a unit ball in $\BC^m$,
there exists $[\sg] \in H^0(Y, F^{n-q})$ such that $\sg|_{X_0} = \sg_0$.
By Lemma \ref{rdh}, we have $\rd_{h_0}(\sg|_{X_0}) = 0$.

(2) 
We take $[\tau] \in H^0(Y, \CK^{n-q})$ such that $\tau|_{X_0} = \tau_0$.
We have $L_f^q([\tau]) = 0$.
Recall the definition of $L_f^q = (\ga^q)^{-1} \circ L^q \circ \gb^0$,
where $\gb^0([\tau]) = \tau \wed dt$,
and $(\ga^q)^{-1}$ is an isomorphism.
Then we have $L^q \circ \gb^0([\tau]) = 0$ in $H^q(X, \Om_X^{n+m}(E))$,
namely $\w^q \wed \tau \wed dt = \rdb a$ for some 
$a \in A^{n+m,q-1}(X, E)$.
By a bidegree reason, $a$ can be written as $a = b \wed dt$ for 
some $b \in A^{n,q-1}(X, E)$. 
Then $(\w^q \wed \tau - \rdb b) \wed dt = 0$.
By restricting on $X_0$, we have $\w_0^q \wed \tau_0 - \rdb b_0 = 0$,
where $b_0 = b|_{X_0}$.

(3)
By (2), we have
$\int_{X_0} \w_0^q \wed \sg_0 \wed h_0 \ol{\tau_0}
= \int_{X_0} \sg_0 \wed h_0 \ol{\rdb b_0}$.
Since $\rd(\sg_0 \wed h_0 \ol{b_0})
= (\rd_{h_0} \sg_0) \wed h_0 \ol{b_0} 
	+ (-1)^{n-q} \sg_0 \wed h_0 \ol{\rdb b_0}$,
and since $\rd_{h_0} \sg_0 = 0$ by (1), we have
$\int_{X_0} \sg_0 \wed h_0 \ol{\rdb b_0} 
= (-1)^{n-q} \int_{X_0} \rd(\sg_0 \wed h_0 \ol{\rdb b_0}) = 0$.
\end{proof}

\subsection{Local freeness}

We shall show that the direct image sheaves $R^qf_* \Om_{X/Y}^{n}(E)$ 
are locally free.
This is an immediate consequence of a result of Takegoshi \cite{Tk}.
We start with recalling a general remark.

\begin{lem} \label{base}
Let $X$ and $Y$ be varieties (reduced and irreducible),
$f : X \lra Y$ be a proper surjective morphism,
and let $\CE$ be a coherent sheaf on $X$ which is flat over $Y$.
Assume that the natural map
$\vph^q(y) : R^qf_*\CE \ot \BC(y) \lra H^q(X_y, \CE_y)$
is surjective for any $y \in Y$ and any $q \geq 0$,
where $X_y$ is the fiber over $y$, and $\CE_y$ is the induced sheaf
(\cite[III.9.4]{Ha}).
Then $R^qf_*\CE$ is locally free for any $q \ge 0$, and
$\vph^q(y) : R^qf_*\CE \ot \BC(y) \lra H^q(X_y, \CE_y)$
is an isomorphism for any $y \in Y$ and any $q \geq 0$.
\end{lem}

\begin{proof}
By \cite[III.12.11(a)]{Ha} (cohomology and base change),
the surjectivity of $\vph^q(y)$ implies that it is an isomorphism.
By \cite[III.12.11(b)]{Ha},
the local freeness of $R^qf_*\CE$ in a neighborhood of $y \in Y$
follows from the surjectivities of $\vph^q(y)$ and of $\vph^{q-1}(y)$.

We can find the corresponding results in the category of complex spaces, 
for example \cite[III.3.4, III.3.7]{BS}.
\end{proof}

\begin{lem} (cf.\ \cite[6.8]{Tk}) \label{free}
Let $f : (X, \w_f) \lra Y$ and $(E, h)$ be as in \S \ref{setup}.I.
Then

(1) the natural restriction map
$R^qf_*\Om_{X/Y}^n(E) \lra H^q(X_y, \Om_{X_y}^n(E_y))$
is surjective for any $y \in Y$ and any $q \ge 0$, and

(2) $R^qf_*\Om_{X/Y}^n(E)$ is locally free for any $q \ge 0$, and
$\vph^q(y) : R^qf_*\Om_{X/Y}^n(E) \ot \BC(y) \lra H^q(X_y, \Om_{X_y}^n(E_y))$
is an isomorphism for any $y \in Y$ and any $q \geq 0$.

\end{lem}

\begin{proof}
(1) Fix $y \in Y$.
Since our assertion is local on $Y$, we may assume that 
$Y$ is a unit ball in $\BC^m$ with coordinates
$t = (t_1, \ldots, t_m)$ centered at $y = \{ t=0 \}$.
We take a trivialization $\CO_Y \cong \Om_Y^m$ given by 
$1 \mapsto dt = dt_1 \wed \ldots \wed dt_m$.
For every $i$ with $1 \le i \le m$, we let 
$Y_i = \{ t_1 = \ldots = t_i = 0 \}$ be a complex sub-manifold of $Y$,
$X_i = f^{-1}(Y_i)$, and let $f_i : X_i \lra Y_i$ be the induced morphism.
We denote by $X_0 = X, Y_0 = Y$ and $f_0 = f$. 
By the injectivity theorem of Takegoshi with $F = \CO_X$ 
in \cite[6.8.i]{Tk}, the sheaf homomorphism
$R^qf_{0*}(f_0^*t_{1}) : 
R^qf_*\Om_{X/Y}^n(E) \ot \Om_Y^m \lra R^qf_*\Om_{X/Y}^n(E) \ot \Om_Y^m$
induced by the product with the holomorphic function $f^*t_1$
is injective for any $q \ge 0$.
Hence the restriction map
$R^qf_{0*}\Om_{X/Y}^n(E) \lra R^qf_{1*}(\Om_{X/Y}^n(E) \ot \CO_{X_1})$
is surjective for any $q \ge 0$.
By the adjunction formula, we have 
$\Om_{X/Y}^n \ot \CO_{X_1} = \Om_{X_1/Y_1}^n$.
Hence inductively, we obtain a surjection
$R^qf_*\Om_{X/Y}^n(E) \lra H^q(X_y, \Om_{X_y}^n(E_y))$.

(2) This follows from (1) and Lemma \ref{base}
\end{proof}


\section{The Hodge metric} \label{hodge}

We shall define a canonical Hermitian metric on $\Rq$,
and compute the metric connection and the curvature.
\S \ref{defhod} will be discussed in the global setting \S \ref{setup}.I,
and the rest of this section will be discussed 
in the localized setting \S \ref{setup}.II.

\subsection{Definition of Hodge metrics} \label{defhod}

We define a canonical Hermitian metric on a vector bundle $\Rq$, 
which we call {\it the Hodge metric}.

\begin{dfn} \label{defhodge}
Let $f :(X, \w_f) \lra Y$ and $(E, h)$ be as in \S \ref{setup}.I,
and let $0 \le q \le n$.
For every point $y \in Y$, we take a local coordinate
$W \cong \{ t \in \BC^m; \ \|t\| < 1 \}$ centered at $y$, 
and a K\"ahler form 
$\w = \w_f + c f^*(\ai \sum dt_j \wed d\ol{t_j})$ on $X_W$
for a real number $c$.
A choice of a K\"ahler form $\w$ gives an injection
$S_\w := S_f^q : R^qf_*\Om_{X/Y}^n(E) \lra f_*\Om_{X/Y}^{n-q}(E)$ over $W$
(Corollary \ref{split}).
Then for every pair of vectors $u_y, v_y \in \Rq_y$, we define
$$
 g(u_y, v_y) 
 = \int_{X_y} (c_{n-q} /q!) 
	(\w_f^q \wed S_\w(u_y) \wed h\ol{S_\w(v_y)})|_{X_y}.
$$
Here the right hand side is the restriction on the image of $S_\w$
of the canonical pairing, say $g$ again, 
on $f_*\Om_{X/Y}^{n-q}(E)_y$ in \S \ref{can}.
\end{dfn}

The injection $S_\w = S_f^q : R^qf_*\Om_{X/Y}^n(E) \lra f_*\Om_{X/Y}^{n-q}(E)$ 
over $W$ may depend on the choices of K\"ahler forms 
in the relative K\"ahler class $\{\w_f\}$, 
however

\begin{lem}
In the situation in Definition \ref{defhodge}, 
the induced metric $g$ on $\Rq|_W$ via 
$S_\w = S_f^q : R^qf_*\Om_{X/Y}^n(E) \lra f_*\Om_{X/Y}^{n-q}(E)$ over $W$
does not depend on the choice of a K\"ahler form $\w$ 
in the relative K\"ahler class $\{\w_f|_{X_W}\} \in H^2(X_W, \BR)$ so that 
$\w|_{X_y} = \w_f|_{X_y}$ for any $y \in W$,
and hence $g$ defines a global Hermitian metric on $\Rq$ over $Y$ 
by varying $y \in Y$.
\end{lem}

\begin{proof}
It is enough to check it in case $u_y = v_y$.
We may also assume that $Y = W \subset \BC^m$.
We take two K\"ahler forms $\w_1$ and $\w_2$ on $X$, which relate
$\w_1 - \w_2 = f^*\ai\rd\rdb \psi$ for some $\psi \in A^0(Y, \BR)$.

(i)
We need to recall the definition of $S_f^q$.
Let $u \in H^0(Y, \Rq) \cong H^q(X, \Om_{X/Y}^{n}(E))$ 
be an extension of $u_y$.
With respect to $\w_i$, we denote by $*_i$ the Hodge $*$-operator,
by $\CH^{n+m,q}(X, \w_i, E, \Phi)$ the space of harmonic forms
in Theorem \ref{Tk4.3} with $\Phi = f^*\|t\|^2$, by
$
\iota_i : \CH^{n+m,q}(X, \w_i, E, \Phi) \isom H^q(X, \Om_X^{n+m}(E))
$
the isomorphism in Theorem \ref{Tk5.2}\,(1),
and by $S_i$ the injection 
$S_f^q : R^qf_*\Om_{X/Y}^n(E) \lra f_*\Om_{X/Y}^{n-q}(E)$.
We have isomorphisms
$$
	\CH_i = \iota_i^{-1} \circ \ga^q : 
	H^q(X, \Om_{X/Y}^n(E)) \isom \CH^{n+m,q}(X, \w_i, E, \Phi)
$$
(see the argument before Corollary \ref{relTk5.2}).
Then we have $*_i\CH_i(u) \in H^0(X, \Om_X^{n+m-q}(E))$,
and $*_i\CH_i(u) = \sg_i \wed dt$ for some 
$[\sg_i] \in H^0(X, \Om_{X/Y}^{n-q}(E))$.
Namely $S_i(u) = [\sg_i]$.
In this setting, our lemma is reduced to show that 
$$
	f_*(\w_1^q \wed \sg_1 \wed h\ol{\sg_1})
	= f_*(\w_2^q \wed \sg_2 \wed h\ol{\sg_2}).	
$$
This is reduced to show 
$\int_X (f^*\gb) \w_1^q \wed \sg_1 \wed dt \wed h \ol{\sg_1 \wed dt}
= \int_X (f^*\gb) \w_2^q \wed \sg_2 \wed dt \wed h \ol{\sg_2 \wed dt}$
for any $\gb \in A^0(Y, \BC)$ with compact support.
We take such a $\gb \in A^0(Y, \BC)$.

(ii)
Since the Dolbeault cohomology classes of $\CH_1(u)$ and
$\CH_2(u)$ are the same, there exists $a \in A^{n+m,q-1}(X, E)$
such that $\CH_1(u) - \CH_2(u) = (c_{n+m-q}/q!) \rdb a$.
Recalling Corollary \ref{relation}\,(2) that 
$\CH_i(u) = (c_{n+m-q}/q!) \w_i^q \wed *_i\CH_i(u)$, we have 
$\w_1^q \wed *_1\CH_1(u) - \w_2^q \wed *_2\CH_2(u) = \rdb a$, and 
hence $\w_1^q \wed \sg_1 \wed dt - \w_2^q \wed \sg_2 \wed dt = \rdb a$.

By a degree reason in the base variables,
we have $f^* (\rdb\gb) \wed d\ol t = 0$. 
Hence
$\rdb((f^*\gb) a \wed h \ol{\sg_1 \wed dt}) 
= (f^*\gb) \rdb a \wed h \ol{\sg_1 \wed dt}
+ (-1)^{n+m+q-1} (f^*\gb) a \wed h \ol{\rd_h(\sg_1 \wed dt)}$.
We also have $\rd_h(\sg_1 \wed dt) = \rd_h *_1\CH_1(u) = 0$ 
by Theorem \ref{Tk4.3}\,(2).
Hence $\int_X (f^*\gb) \rdb a \wed h \ol{\sg_1 \wed dt}
= \int_X \rdb((f^*\gb) a \wed h \ol{\sg_1 \wed dt}) = 0$
by the Stokes theorem.
Then the relation
$\w_1^q \wed \sg_1 \wed dt = \w_2^q \wed \sg_2 \wed dt + \rdb a$
implies that
$\int_X (f^*\gb) \w_1^q \wed \sg_1 \wed dt \wed h \ol{\sg_1 \wed dt}
= \int_X (f^*\gb) \w_2^q \wed \sg_2 \wed dt \wed h \ol{\sg_1 \wed dt}$.

(iii)
Now we use $\w_1 - \w_2 = f^*\ai\rd\rdb \psi$. 
This leads
$\int_X (f^*\gb) \w_2^q \wed \sg_2 \wed dt \wed h \ol{\sg_1 \wed dt}
= \int_X (f^*\gb) \w_1^q \wed \sg_2 \wed dt \wed h \ol{\sg_1 \wed dt}$.
The last integral equals to
$\int_X (f^*\gb) \sg_2 \wed dt \wed h \ol{\w_1^q \wed \sg_1 \wed dt}
= \int_X (f^*\gb) \sg_2 \wed dt \wed h \ol{\w_2^q \wed \sg_2 \wed dt}
+ \int_X (f^*\gb) \sg_2 \wed dt \wed h \ol{\rdb a}$.
By a similar manner as above, mainly because of 
$\rd_h(\sg_2 \wed dt) = \rd_h *_2\CH_2(u) = 0$, we can see $
\int_X (f^*\gb) \sg_2 \wed dt \wed h \ol{\rdb a} = 0$.
We finally obtain
$\int_X (f^*\gb) \w_1^q \wed \sg_1 \wed dt \wed h \ol{\sg_1 \wed dt}
= \int_X (f^*\gb) \w_2^q \wed \sg_2 \wed dt \wed h \ol{\sg_2 \wed dt}$.
\end{proof}

At this point, we have the so-called the metric connection 
(or the Chern connection) $D_g$ of the Hermitian vector bundle $(\Rq, g)$, 
and the curvature $\Th_g = D_g^2$.
Since the curvature property in Theorem \ref{MT} is a local question 
on the base $Y$, it is enough to consider in the following setting:

Let $f :(X, \w) \lra Y \subset \BC^m$ and $(E, h)$ 
be as in \S \ref{setup}.II, and let $0 \le q \le n$.
We denote by $F = F^{n-q}$ the image of 
$S_f^q : R^qf_*\Om_{X/Y}^{n}(E) \lra f_*\Om_{X/Y}^{n-q}(E)$
with respect to $\w$.
Since, by definition, the canonical pairing $g$ on $f_*\Om_{X/Y}^{n-q}(E)$
gives our $(\Rq, g)$, we say a sub-bundle 
$$
	(F, g) \subset (f_*\Om_{X/Y}^{n-q}(E), g).
$$

\subsection{The metric connection}


We shall construct the metric connection $D_g$ of $(F, g)$.
Recall Lemma \ref{rdh} that $\rd_h \sg = \sum_j \mu^j \wed dt_j$ with some 
$\mu^j \in A^{n-q,0}(X, E)$ for $[\sg] \in H^0(Y, F)$.
Since $A^0(Y, F) = A^0(Y, \BC) \ot H^0(Y, F)$ as $A^0(Y, \BC)$-module,
this formula holds for $[\sg] \in A^0(Y, F)$, too.
We consider the fiberwise orthogonal projection
$P_y : A^{n-q,0}(X_y, E_y) \lra F_y$
given by $u_y \mapsto \sum_{j=1}^\ell g_y(u_y, \sg_{jy})\sg_{jy}$,
where $\sg_{1y}, \ldots, \sg_{\ell y} \in F_y$ is a basis of $F_y$.
Since $F$ is locally free, the family $\{P_y\}_{y \in Y}$ induces a map
$$
	P : A^{n-q,0}(X, E) \lra 
\{ u \in A^{n-q,0}(X, E); \ u|_{X_y} \in F_y \text{ for any } y \in Y \} 
$$
Then for $[\sg] \in A^0(Y, F)$ with $\rd_h \sg = \sum_j \mu^j \wed dt_j$, 
we define
$$
	\rd_g [\sg] = \sum [P(\mu^j)] dt_j \in A^{1,0}(Y, F).
$$

\begin{lem}
The class $[P(\mu^j)]$ is well-defined for $[\sg] \in A^0(Y, F)$.
\end{lem}

\begin{proof}
(1)
We shall show that $\mu^j|_{X_y}$ are perpendicular to
$H^{0}(X_y, \Om_{X_y}^{n-q}(E_y))$ under the condition $[\sg] = [0]$, namely
$\sg|_{X_y} = 0$ for any $y \in Y$.
We write as $\sg = \sum \sg_j \wed dt_j$ with some 
$\sg_j \in A^{n-q-1,0}(X, E)$.
We note that we can take $\mu^j = \rd_h \sg_j$. 
We take any $s \in H^{0}(X_y, \Om_{X_y}^{n-q}(E_y))$.
Then $\rd(\w_y^q \wed  \sg_j|_{X_y} \wed h_y \ol{s})
= \w_y^q \wed \rd_{h_y}(\sg_j|_{X_y}) \wed h_y \ol{s}
	+(-1)^{n-q} \w_y^q \wed \sg_j|_{X_y} \wed h_y \ol{\rdb s}$.
Because of $\rdb s = 0$, we have 
$g_y ((\rd_{h}\sg_j)|_{X_y}, s)
= (c_{n-q}/q!) \int_{X_y} \w_y^q \wed \rd_{h_y}(\sg_j|_{X_y}) \wed h_y \ol{s}
= (c_{n-q}/q!) \int_{X_y} \rd(\w_y^q \wed \sg_j|_{X_y} \wed h_y \ol{s}) = 0$.

(2)
The above (1) is enough to show that $[P(\mu^j)]$ is well-defined.
But in fact, (1) said slightly more.
\end{proof}

\begin{lem} 
The sum $D_g := \rd_g + \rdb$ is the metric connection of 
the Hermitian vector bundle $(F, g)$.
\end{lem}

\begin{proof}
It is not difficult to see that it is a connection.
Let us check the compatibility with the metric $g$.
Let $[\sg], [\tau] \in H^0(Y, F)$, and write 
$\rd_h \tau = \sum_j \mu^j(\tau) \wed dt_j$.
Then 
$\rdb g([\sg], [\tau]) 
= (-1)^{n-q} f_*((c_{n-q}/q!) \w^q \wed \sg \wed h \ol{\rd_h \tau}) 
= \sum_j f_*((c_{n-q}/q!) \w^q \wed \sg \wed h \ol{\mu^j(\tau)}) d\ol{t_j}$.
Since $\sg|_{X_y} \in F_y$, the last term becomes
$\sum_j f_*((c_{n-q}/q!) \w^q \wed \sg \wed h \ol{P(\mu^j(\tau))}) d\ol{t_j}$,
and it is  
$\sum_j g([\sg], [P(\mu^j(\tau))]) d\ol{t_j}$.
In the notation of \S \ref{hvb}, we can write as 
$g([\sg], [\tau]) = [\sg] \wed g\ol{[\tau]}$
and $\rdb g([\sg], [\tau]) = [\sg] \wed g \ol{\rd_g[\tau]}$.
%
\end{proof}

\subsection{Curvature formula}

We describe the Nakano semi-positivity of 
a Hermitian holomorphic vector bundle.
Since it is a local property, we will discuss on a local chart.
Let $Y \subset \BC^m$ be a unit ball centered at $0$
with coordinates $t = (t_1, \ldots, t_m)$,
and let $F = Y \times \BC^\ell$ be a trivial vector bundle
with a non-trivial Hermitian metric $g$.
(This $(F, g)$ may not necessarily be our original bundle.) \ 
We write $\Th_g = \sum \Th_{jk} dt_j \wed d\ol{t_k}$
with $\Th_{jk} \in \End(Y, F)$.
Then $(F, g)$ is Nakano semi-positive at $t = 0$, if and only if
for any tensor $s = \sum_{j=1}^m \rd/\rd t_j \ot \sg^0_{j}
\in (T_Y \ot F)_{0}$, we have
$\Th_g(s) = \sum_{j,k} g_0(\Th_{jk}\sg^0_{j}, \sg^0_{k}) \ge 0$.
Moreover the last quantity can be obtained another way from
local sections.
If $\sg, \tau \in H^0(Y, F)$, we have 
$\frac{\rd^2}{\rd t_j \rd \ol{t_k}} g(\sg, \tau)
= g((\rd_g \sg)^j, (\rd_g \tau)^k) - g(\Th_{jk}\sg, \tau)$,
where $\rd_g \sg = \sum_j (\rd_g \sg)^j dt_j \in A^{1,0}(Y, F)$ and so on.
Hence if $\sg$ and $\tau$ are {\it normal} at $0$ with respect to $g$
(i.e., $\rd_g \sg = \rd_g \tau = 0$ at $0$), we have
$(\frac{\rd^2}{\rd t_j \rd \ol{t_k}} g(\sg, \tau))|_{t=0}
= - g_0(\Th_{jk} \sg|_{t=0}, \tau|_{t=0})$.

\begin{notation}
(1) Let $V$ be a continuous $(m,m)$-form on $Y \subset \BC^m$.
Then we can write $V = v(t) dV_t$ with a unique continuous function $v$ on $Y$,
and we define $V_{t=0} := v(0)$.

(2) Associated to $m$-ple $\sg_1, \ldots, \sg_m \in H^0(Y, F)$,
we let 
$$
	T(\sg) = \sum_{j,k} g(\sg_j, \sg_k) \what{dt_j \wed d\ol{t_k}}
		\ \in A^{m-1,m-1}(Y, \BC).
$$
\end{notation}

In case all $\sg_j$ are normal at $t = 0$, we have
$\levi T(\sg)_{t=0} = - \sum_{j,k} g_0(\Th_{jk}\sg_j|_{x_0}, \sg_k|_{x_0})$.
Hence we have

\begin{lem} \cite[\S 2]{B} \label{cri}
A Hermitian vector bundle $(F, g)$ on an open subset $Y \subset \BC^m$ 
is Nakano semi-positive at $t=0$, if
for any $m$-ple vectors $\sg^0_{1}, \ldots, \sg^0_{m} \in F_0$, 
there exist extensions $\sg_j \in H^0(Y, F)$ of $\sg^0_{j}$,
all of which are normal at $t = 0$ and satisfy $\levi T(\sg)_{t=0} \leq 0$. 
\end{lem}

We go back to our original situation.
We prepair the following notations.

\begin{notation} \label{notation}
Let $f : (X, \w) \lra Y \subset \BC^m$ and $(E, h)$ be as in \S \ref{setup}.II.
Let $\sg_1, \ldots, \sg_m \in A^{n-q,0}(X, E)$ such that 
$[\sg_j] \in H^0(Y, F)$ for all $j$.

(1) We set 
$$
	\hsg = \sum \sg_j \wed \what{dt_j} \in A^{n-q+m-1,0}(X, E).
$$
Then
$$
	T([\sg]) = \sum_{j,k} g([\sg_j], [\sg_k]) \what{dt_j \wed d\ol{t_k}}
		 = f_*((c_N/q!) \w^q \wed \hsg \wed h \ol\hsg ).
$$
Here $N = n-q+m-1$.

(2) We write $\rd_h \sg_j = \sum_k \mu_j^k \wed dt_k$.
Then 
$$
	\rd_h \hsg = \sum_j \mu_j^j \wed dt =: \mu \wed dt
$$
with $\mu \in A^{n-q,0}(X, E)$, or rather $[\mu] \in A^{n-q,0}(X/Y,E)$.

(3) We write $\rdb \sg_j = \sum_k \eta_j^k \wed dt_k$.
Then 
$$
	\rdb \hsg = \sum_j \eta_j^j \wed dt =: \eta \wed dt
$$
with $\eta \in A^{n-q-1,1}(X, E)$, or rather $[\eta] \in A^{n-q-1,1}(X/Y,E)$.
\qed
\end{notation}

\begin{lem} \label{f1} 
(cf.\ \cite[(4.4)]{B}) \ 
In Notation \ref{notation}, one has
\begin{eqnarray*}
-\ai\rd\rdb T([\sg])_{t=0}
& = & f_* ((c_N/q!) \w^q \wed \ai \Th_h \wed  \hsg \wed h \ol\hsg)_{t=0}\\
& &\ \  - \int_{X_0} (c_{n-q}/q!) (\w^q \wed  \mu \wed h \ol\mu)|_{X_0} 
	- \int_{X_0} (c_{n-q}/q!) (\w^q \wed  \eta \wed h \ol\eta)|_{X_0}.
\end{eqnarray*}
\end{lem}

\begin{rem}
The first term comes from the curvature of $E$, and contributes positively.

The second term is $-\|\mu|_{X_0}\|_{h_0}^2$, and it can be seen as 
the ``second fundamental form'' of 
$F \subset \bigcup_{t \in Y} A^{n-q,0}(X_y,E_y)$ at $t=0$.
This negative contribution will be eliminated by a careful choice
of forms $\sg_j$, in \S \ref{exist}.

The third term is not a definite form.
In general one can write $\eta|_{X_0}$ as a sum
$\eta|_{X_0} = \eta_0' + \w_0 \wed \eta_0''$ 
for primitive forms $\eta_0'$ and $\eta_0''$ on $X_0$, and then  
$- \int_{X_0} (c_{n-q}/q!) (\w^q \wed  \eta \wed h \ol\eta)|_{X_0}
= \|\eta_0'\|_{h_0}^2 - \|\eta_0''\|_{h_0}^2$.
In \S \ref{exist}, we will show that we can take $\sg_j$ so that
all $\eta_j^k|_{X_0}$ and hence $\eta|_{X_0}$ are primitive on $X_0$.
In that case, the third term is 
$- \int_{X_0} (c_{n-q}/q!) (\w^q \wed  \eta \wed h \ol\eta)|_{X_0}
= \|\eta|_{X_0}\|_{h_0}^2 \geq 0$. 
We should read the Kodaira-Spencer class contributes positively.
\end{rem}

\begin{proof}[Proof of Lemma \ref{f1}]
The proof will be done by direct computations.
We first note that 
$f_*(\w^q \wed \rdb \hsg \wed h \ol\hsg)
= f_*(\w^q \wed \eta \wed dt \wed h \ol\hsg) = 0$
as an $(m-1,m)$-current on $Y$, because it contains $dt$.
By the same reason, we have
$f_*(\w^q \wed \hsg \wed h \ol{\rdb \hsg})=0$,
and hence, by taking $\rdb$, we have
$f_*(\w^q \wed \hsg \wed h \ol{\rd_h \rdb \hsg})
= - (-1)^N f_*(\w^q \wed \rdb \hsg \wed h \ol{\rdb \hsg})$.
Then we have
$\rdb f_*(\w^q \wed \hsg \wed h \ol\hsg) 
= (-1)^N f_*(\w^q \wed \hsg \wed h \ol{\rd_h \hsg})$, and then
$$
\rd \rdb f_*(\w^q \wed \hsg \wed h \ol\hsg) 
= (-1)^N f_*(\w^q \wed \rd_h \hsg \wed h \ol{\rd_h \hsg})
 + f_*(\w^q \wed \hsg \wed h \ol{\rdb \rd_h \hsg}).
$$
Since $\rd_h\rdb + \rdb \rd_h = e(\Th_h)$, we have
$f_*(\w^q \wed \hsg \wed h \ol{\rdb \rd_h \hsg})
=  f_*(\w^q \wed \hsg \wed h \ol{\Th_h \wed \hsg})
 - f_*(\w^q \wed \hsg \wed h \ol{\rd_h \rdb \hsg})$.
Using
$f_*(\w^q \wed \hsg \wed h \ol{\rd_h \rdb \hsg})
= - (-1)^N f_*(\w^q \wed \rdb \hsg \wed h \ol{\rdb \hsg})$,
we can write
\begin{equation*} 
\begin{aligned}
\rd \rdb f_*(\w^q \wed \hsg \wed h \ol\hsg)
 = 
 - f_*(\w^q \wed \Th_h \wed \hsg \wed h \ol\hsg) 
& +(-1)^{N+(n-q)m} f_*(\w^q \wed \mu \wed h \ol\mu \wed dt \wed d\ol t) \\
& +(-1)^{N+(n-q)m} f_*(\w^q \wed \eta \wed h \ol\eta \wed dt \wed d\ol t).
\end{aligned}
\end{equation*} 
Here we mind that $\ai\Th_h$ is real.
Hence
$-\ai\rd\rdb f_*((c_N/q!) \w^q \wed \hsg \wed h \ol\hsg)$ 
is
\begin{equation*} 
\begin{aligned}
f_*((c_N/q!) \w^q \wed \ai \Th_h \wed \hsg \wed h \ol\hsg) 
& - f_*((c_{n-q}/q!) \w^q \wed \mu \wed h \ol\mu \wed c_m dt \wed d\ol t) \\
& - f_*((c_{n-q}/q!) \w^q \wed \eta \wed h \ol\eta \wed c_m dt \wed d\ol t).
\end{aligned}
\end{equation*} 
By taking their values at $t = 0$, we have our assertion. 
\end{proof}


\section{Normal and ``primitive'' sections, 
and the proof of Theorem \ref{MT}}

Let $f : (X, \w) \lra Y \subset \BC^m$ and $(E, h)$ as in \S \ref{setup}.II,
and keep the notations in \S \ref{hodge}.

\subsection{Effect of normality}

We control $\rd_h \sg$ at one point for $[\sg] \in H^0(Y, F)$,
when it is normal at $t = 0$.
Recall $\rd_h \sg = \sum \mu^j \wed dt_j$ with some $\mu^j \in A^{n-q,0}(X,E)$.
To go further, we need a genericity condition over the base $Y$.
We will assume that the function 
$y \mapsto \dim H^0(X_y, \Om_{X_y}^{n-q}(E_y))$ is constant around $t = 0$.
This assumption implies that $f_*\Om_{X/Y}^{n-q}(E)$ is locally free around
$t=0$, and that the fiber $f_*\Om_{X/Y}^{n-q}(E)_y$ coincides with
$H^0(X_y, \Om_{X_y}^{n-q}(E_y))$ around $t=0$ 
(\cite[10.5.5]{GR}, \cite[III \S 12]{Ha}).
In case $q = 0$, i.e., the case in \cite[\S 4]{B},
this assumption holds true thanks to Ohsawa-Takegoshi type $L^2$-extension
theorem \cite{OT}\,\cite{O}\,\cite{Ma}.
Recall Corollary \ref{split} that $f_*\Om_{X/Y}^{n-q}(E) = F \oplus \CK$, 
where $\CK = \CK^{n-q} = \Ker L_f^q$.

%

\begin{lem} \label{s-normal}
Assume that the function $y \mapsto \dim H^0(X_y, \Om_{X_y}^{n-q}(E_y))$
is constant around $t = 0$.
Let $[\sg] \in H^0(Y, F)$ with $\rd_h \sg = \sum \mu^j \wed dt_j$, 
and suppose $\rd_g[\sg]=0$ at $t = 0$.
Then all $\mu^j|_{X_0}$ are perpendicular to $H^0(X_0, \Om_{X_0}^{n-q}(E_0))$.
\end{lem}

\begin{proof}
We will use notations in \S \ref{abs} for $(X_0, \w_0)$ and $(E_0, h_0)$.
Let $(\ ,\ )_{h_0}$ be the inner product of $A^{n-q,0}(X_0, E_0)$ 
in terms of the metrics $\w_0$ and $h_0$ on $X_0$.
We have $H^0(X_0, \Om_{X_0}^{n-q}(E_0)) = F_0 \oplus \CK_0$,
which is an orthogonal direct sum by our assumption and by 
Lemma \ref{ortho}\,(3).
We fix $j$.
Let $\mu^j|_{X_0} = \tau_0 + a_0 \in A^{n-q,0}(X_0, E_0)$ be
the Hodge decomposition of forms so that 
$\tau_0 \in H^0(X_0, \Om_{X_0}^{n-q}(E_0))$ and 
$a_0 \in \vth_{h_0}A^{n-q,1}(X_0, E_0)$.
We would like to show that $\tau_0 = 0$.

Since $\rd_g[\sg]=0$ at $t = 0$, namely all $\mu^k|_{X_0}$ are perpendicular
to $F_0$, it follows that $\tau_0 \in \CK_0$ by Lemma \ref{ortho}\,(3). 
Then $\w_0^q \wed \tau_0 = \rdb b_0$ for some $b_0 \in A^{n, q-1}(X_0, E_0)$ 
by Lemma \ref{ortho}\,(2).
Combining with Lemma \ref{rdh} that $\rd_{h_0}(\mu^j|_{X_0}) = 0$, we have
$\int_{X_0} \w_0^q \wed \mu^j|_{X_0} \wed h_0\ol{\tau_0} = 0$ 
by integration by parts as in Lemma \ref{ortho}\,(3).
Then
$\|\tau_0\|_{h_0}^2 = (\tau_0+a_0,\tau_0)_{h_0} 
= \int_{X_0} (c_{n-q}/q!) \w_0^q \wed \mu^j|_{X_0} \wed h_0\ol{\tau_0} = 0$,
and hence $\tau_0 = 0$. 
\end{proof}

\begin{lem} \label{rd-exact}
Let $[\sg] \in H^0(Y, F)$ with $\rd_h \sg = \sum \mu^j \wed dt_j$,
and suppose that $\mu^j|_{X_0}$ is perpendicular to 
$H^0(X_0, \Om_{X_0}^{n-q}(E_0))$ for some $j$.
Then there exists $\xi_0^j \in A^{n-q-1,0}(X_0, E_0)$ such that
$\rd_{h_0} \xi_0^j = \mu^j|_{X_0}$ and that $\rdb \xi_0^j$ is primitive.
\end{lem}

\begin{proof}
We will use notations in \S \ref{abs} for $(X_0, \w_0)$ and $(E_0, h_0)$.
Recall Theorem \ref{Tk4.3} with $\dim Y = 0$ that 
the Hodge $*$-operator yields an injective homomorphism
$*_0 : \CH^{n,q}(X_0, E_0) \lra H^0(X_0, \Om_{X_0}^{n-q}(E_0))$.

We consider $u := \w_0^q \wed \mu^j|_{X_0} \in A^{n,q}(X_0, E_0)$,
and recall $(c_{n-q}/q!) *_0 u = \mu^j|_{X_0}$.
Let $u = a + \rdb b + \vth_{h_0} c$ be the Hodge decomposition
of forms so that $a \in \CH^{n,q}(X_0, E_0), 
b$ is $\vth_{h_0}$-exact, and that $c$ is $\rdb$-exact.

We first show that $\rdb b = 0$.
Using integration by parts and by Lemma \ref{rdh}\,(3), we have
$\int_{X_0} \rdb b \wed h_0 \ol{\mu^j}|_{X_0} = 0$.
Since $\| \rdb b\|_{h_0}^2 
= (\rdb b, u)_{h_0}
= \int_{X_0} \rdb b \wed \ol{*}_{h_0} u$,
and since the last term is
$c_{n-q}q! \int_{X_0} \rdb b \wed h_0 \ol{\mu^j}|_{X_0} = 0$,
we have $\rdb b = 0$.

We next show that $a = 0$.
Recall in general, $(v, w)_{h_0} = (*_0 v, *_0 w)_{h_0}$ holds
for $v, w \in A^{p,q}(X_0, E_0)$ (\cite[1.2.20]{Huy}).
Since $\vth_{h_0} c \in (\CH^{n,q}(X_0, E_0))^\perp$
the orthogonal complement in $A^{n,q}(X_0, E_0)$,
we have $*_0 (\vth_{h_0} c) \in (*_0 \CH^{n,q}(X_0, E_0))^\perp$.
We also have $*_0 u = c_{n-q}^{-1}q! \mu^j|_{X_0}
\in H^0(X_0, \Om_{X_0}^{n-q}(E_0))^\perp 
\subset (*_0 \CH^{n,q}(X_0, E_0))^\perp$.
On the other hand $*_0a \in *_0 \CH^{n,q}(X_0, E_0)$,
hence the both sides of $*_0a = *_0u -*_0(\vth_{h_0} c)$ have to be $0$.

Now we had $u = \vth_{h_0} c$ for a $\rdb$-exact form 
$c \in A^{n,q+1}(X_0, E_0)$.
By the Lefschetz isomorphism on forms (\cite[1.2.30]{Huy}), there exists 
$\xi \in A^{n-q-1,0}(X_0, E_0)$ such that $\w_0^{q+1} \wed \xi = c$.
We have $\w_0^{q+1} \wed \rdb \xi = \rdb c = 0$, namely 
$\rdb \xi$ is primitive.
We also have $*_0c = *_0(\w_0^{q+1} \wed \xi) 
= c_{n-q-1}^{-1}(q+1)! \xi$ by Lemma \ref{primitive}.
Then 
$\mu^j|_{X_0} = (c_{n-q}/q!) *_0 u 
= (c_{n-q}/q!) (-*_0 \circ *_0 \rd_{h_0}(*_0 c))
= -(-1)^{n-q} (c_{n-q}/q!) \rd_{h_0} (c_{n-q-1}^{-1}(q+1)! \xi)
= -\ai (q+1) \rd_{h_0} \xi$.
We finally take $\xi_0 = -\ai (q+1)\xi$.
\end{proof}

\subsection{Existence of strongly normal and ``primitive'' sections} 
\label{exist}

Here we state a key result for the curvature estimate of our Hodge metric,
as a consequence of Lemma \ref{rdb} and Lemma \ref{rd-exact}.
The part (I) of Proposition \ref{strongNP} below in fact holds 
not only for $F$, but also for any locally free subsheaf 
of $f_*\Om_{X/Y}^{n-q}(E)$.
The property (3) (respectively, (4)) below will be referred as 
``primitive'' (respectively, strongly normal) at $t=0$.

\begin{prop} \label{strongNP} (cf.\ \cite[Proposition 4.2]{B}) \ 
Let $\sg^0 \in F_0$ be a vector at $t=0$.

(I) 
Then, there exists $\sg \in A^{n-q, 0}(X, E)$ such that 
$[\sg] \in H^0(Y, F)$ with the following three properties:\
(1)  $\sg|_{X_0} = \sg^0$,

(2) $\rd_g[\sg] = 0$ at $t=0$, 

(3) $\eta^j|_{X_0} \wed \w_0^{q+1} = 0$ for any $j$,
where $\rdb \sg = \sum \eta^j \wed dt_j$.

(II)
If the function $y \mapsto \dim H^0(X_y, \Om_{X_y}^{n-q}(E_y))$
is constant around $t = 0$, one can take $\sg$ in (I) with
the following additional fourth property:\

(4) $\mu^j|_{X_0} = 0$ for any $j$,
where $\rd_h \sg = \sum \mu^j \wed dt_j$.
\end{prop}

\begin{proof}
(I)
A local extension as in (1) and (2) is possible for any Hermitian
vector bundle.
Hence we start with a local extension $[\sg] \in H^0(Y,F)$ 
satisfying (1) and (2).
We write $\rdb \sg = \sum \eta^j \wed dt_j$.
By Lemma \ref{rdb}\,(4), we have $(\eta^j \wed \w^{q+1})|_{X_0} = \rdb a^j_0$
for some $a^j_0 \in A^{n,q+1}(X_0,E_0)$.
By the Lefschetz isomorphism on forms (\cite[1.2.30]{Huy}), we can write
$a^j_0 = b^j_0 \wed \w_0^{q+1}$ for some $b^j_0 \in A^{n-q-1,0}(X_0,E_0)$.
We take smooth extensions $b^j \in A^{n-q-1,0}(X, E)$ so that
$b^j|_{X_0} = b^j_0$, and we let 
$\wtil \sg = \sg - \sum b^j \wed dt_j \in A^{n-q,0}(X, E)$.
We check $[\wtil \sg]$ is what we are looking for.
Since $[\wtil \sg] = [\sg]$ in $A^{n-q,0}(X/Y, E)$,
we see $[\wtil \sg] \in H^0(Y,F)$, and (1) and (2) for $\wtil \sg$.
Moreover $\rdb \wtil \sg = \sum (\eta^j - \rdb b^j) \wed dt_j$,
and 
$(\eta^j - \rdb b^j)|_{X_0}  \wed \w_0^{q+1}  
= \rdb a^j_0 - \rdb (b_0^j \wed \w_0^{q+1})=0$.
Hence we have $\eta^j(\wtil \sg)|_{X_0} \wed \w_0^{q+1} = 0$,
i.e., (3) for $\wtil \sg$.

(II)
We assume that the function $y \mapsto \dim H^0(X_y, \Om_{X_y}^{n-q}(E_y))$
is constant around $t = 0$.
We take $\sg \in A^{n-q,0}(X, E)$ which satisfies all three properties in (I).
We write $\rd_h \sg = \sum \mu^j \wed dt_j$ and 
$\rdb \sg = \sum \eta^j \wed dt_j$.
By Lemma \ref{s-normal} and \ref{rd-exact}, for every $j$,
there exists $\xi_0^j \in A^{n-q-1,0}(X_0, E_0)$ such that
$\rd_{h_0} \xi_0^j = \mu^j|_{X_0}$ and that $\rdb \xi_0^j$ is primitive.
We take $\xi^j \in A^{n-q-1,0}(X, E)$ such that
$\xi^j|_{X_0} = \xi_0^j$ for every $j$.
We consider $\wtil \sg = \sg - \sum_j \xi^j \wed dt_j$.
Since $[\wtil \sg] = [\sg]$ in $A^{n-q,0}(X, E)$, 
we see $[\wtil \sg] \in H^0(Y, F)$, and (1) and (2) for $\wtil \sg$.
We have 
$\rd_h \wtil \sg = \sum (\mu^j -  \rd_h \xi^j) \wed dt_j$ and
$\rdb \wtil \sg = \sum (\eta^j - \rdb \xi^j) \wed dt_j$.
The property (3) for $\wtil \sg$ follows from the primitivity of
$\eta^j|_{X_0}$ for $\sg$ and of $(\rdb \xi^j)|_{X_0} = \rdb \xi_0^j$.
The property (4) for $\wtil \sg$ follows from 
$\mu^j|_{X_0} -  (\rd_h \xi^j)|_{X_0} = 0$.
\end{proof}


\subsection{Nakano semi-positivity}

\begin{prop} \label{negative}
Let $\sg_1, \ldots, \sg_m \in A^{n-q,0}(X, E)$ with
$[\sg_1], \ldots, [\sg_m] \in H^0(Y, F)$ and satisfying the properties
(3) and (4) in Proposition \ref{strongNP}.
Then $\levi T([\sg])_{t=0} \le 0$ in Notation \ref{notation}
for these $\sg_1, \ldots, \sg_m$. 
\end{prop}

\begin{proof}
We will use the notations in \ref{notation}.
By the property (4) in Proposition \ref{strongNP}, we have $\mu|_{X_0} = 0$.
The property (3) in Proposition \ref{strongNP} implies that
$\eta|_{X_0}$ is primitive.
In particular, by using Lemma \ref{primitive},
$- \int_{X_0}(c_{n-q}/q!)(\w^q \wed \eta \wed h \ol \eta)|_{X_0} 
= \| \eta|_{X_0} \|_{h_0}^2$ the square norm with respect to $\w_0$ and $h_0$.
Then the formula in Lemma \ref{f1} is 
$$
-\ai\rd\rdb T([\sg])_{t=0}
= f_* ((c_N/q!) \w^q \wed \ai \Th_h \wed  \hsg \wed h \ol\hsg)_{t=0}
	+ \| \eta|_{X_0} \|_{h_0}^2.
$$
The right hand side is non-negative, since the curvature $\Th_h$ is
Nakano semi-positive.
\end{proof}

\begin{cor}
$(F, g)$ is Nakano semi-positive, and hence so is $R^qf_*\Om_{X/Y}^n(E)$.
\end{cor}

\begin{proof}
Since $g$ is a smooth Hermitian metric of $F$,
to show the Nakano semi-positivity, it is enough to show it
on the complement of an analytic subset of $Y$.
By Grauert (\cite[10.5.4]{GR}\,\cite[III.12.8, 12.9]{Ha}), 
there exists an analytic subset $Z \subset Y$ such that
the function $y \mapsto \dim H^0(X_y, \Om_{X_y}^{n-q}(E_y))$
is constant on $Y \setminus Z$.
We apply the criterion in Lemma \ref{cri} at each point on $Y \setminus Z$.
Then thanks to Proposition \ref{strongNP}, Proposition \ref{negative} 
in fact shows that $g$ is Nakano semi-positive on $Y \setminus Z$.
\end{proof}



\baselineskip=14pt

\vskip 10pt

Christophe MOUROUGANE

\vskip 5pt

Institut de Recherche Math\'ematique de Rennes

Campus de Beaulieu

35042 Rennes cedex, France

e-mail: christophe.mourougane@univ-rennes1.fr

\vskip 10pt

Shigeharu TAKAYAMA

\vskip 5pt

Graduate School of Mathematical Sciences

University of Tokyo

3-8-1 Komaba, Tokyo

153-8914, Japan

e-mail: taka@ms.u-tokyo.ac.jp

\end{document}